\documentclass[11pt,a4paper,twoside]{amsart}

\usepackage{amsmath,amssymb,amsfonts,amsthm,mathrsfs}
\usepackage{times,hyperref,color}

\let\mathcal\mathscr

\newtheorem{The}{Theorem}[section]

\newtheorem{Theorem}{Theorem}[section]

\newtheorem{Proposition}[The]{Proposition}
\newtheorem{Conjecture}[The]{Conjecture}
\newtheorem{Lemma}[The]{Lemma}

\theoremstyle{definition}
\newtheorem{Definition}[The]{Definition}

\newtheorem{Remark}[The]{Remark}


\begin{document}


\title{
On the maximum conjecture
}

\author{Masoud Sabzevari}
\address{Department of Mathematics,
 Shahrekord University, 88186-34141 Shahrekord, IRAN and School of
Mathematics, Institute for Research in Fundamental Sciences (IPM), P.
O. Box: 19395-5746, Tehran, IRAN}
\email{sabzevari@math.iut.ac.ir}

\date{\number\year-\number\month-\number\day}

\maketitle

\begin{abstract}
We verify the maximum conjecture on the rigidity of totally nondegenerate model CR manifolds in the following two cases: (i) for all models of CR dimension one (ii) for the so-called full-models, namely those in which their associated symbol algebras are free CR. In particular, we discover that in each arbitrary CR dimension and length $\geq 3$, there exists at least one totally nondegenerate model, enjoying this conjecture. Our proofs rely upon some recent results in the Tanaka theory of transitive prolongation of fundamental algebras.
\end{abstract}

\pagestyle{headings} \markright{On the maximum conjecture}

\section{Preliminaries}
\label{sec-prel}

For an arbitrary smooth real manifold $M$, a subbundle $T^cM\subset TM$ of  even rank is an {\sl almost Cauchy-Riemann (CR {\rm for short}) structure} if it is equipped with a fiber preserving {\sl complex structure} $J:T^cM\rightarrow T^cM$ {\it i.e.} a linear map enjoying $J\circ J=-{\rm id}$. In this case, $M$ is called an {\sl almost CR manifold} of CR dimension $n:=\frac{1}{2}\cdot({\rm dim}\,T^cM)$ and codimension $k:={\rm dim}\, M-2n$. As is known (\cite{BER, 5-cubic}), the {\sl complexified bundle} $\mathbb C\otimes T^cM$ decomposes as:
\[
\mathbb C\otimes T^cM:=T^{1,0}M\oplus T^{0,1}M,
\]
where:
\[
T^{1,0}M:=\Big\{X-i\,J(X): \ X\in T^cM\Big\}
\]
and $T^{0,1}M=\overline{T^{1,0}M}$ are the {\sl holomorphic and anti-holomorphic subbundles} of $\mathbb C\otimes TM$. By definition, $M$ is a {\sl CR manifold} with {\sl CR structure} $T^cM$ if $T^{1,0}M$ is involutive in the sense of Frobenius. Such a CR manifold is {\sl generic} in $\mathbb C^{n+k}$ if it can be represented as the graph of some $k$ functionally independent defining equations.

Let $p\in M$ be an arbitrary point. The {\sl Levi form} ${\bf L}_p$ of $M$ at this point is the Hermitian form on $T^{1,0}_pM$ with values in $\mathbb C\otimes \frac{T_pM}{T^c_pM}$ defined as:
\[
{\bf L}_p(L_1(p), L_2(p)):=i\,[L_1, \overline L_2](p) \ \ \ \ \ \ {\rm mod} \ \ \mathbb C\otimes T^c_pM,
\]
for each two local vector fields $L_1$ and $L_2$ of $T^{1,0}M$ near $p$.
A CR manifold is {\sl nondegenerate} whenever its associated Levi form is nondegenerate at each point.

The interesting class of nondegenerate CR manifolds is investigated widely during the recent century ({\it see e.g.} \cite{BER, Chern-Moser, Merker-Porten-2006}). Among them, we have {\it totally nondegenerate} CR manifolds which, in some sense, are the {\it most nondegenerate} ones. In order to define total nondegeneracy, consider $ M\subset\mathbb C^{n+k}$ as a real analytic generic submanifold of CR dimension $n$, codimension $k$ and hence of real dimension $2n+k$. With the customary abuse of notation, in this paper we use the same symbol $T^cM$ for the collection of all local vector fields associated to this distribution. Set $D_1:=T^c M$ and define successively $D_j=D_{j-1}+[D_1,D_{j-1}]$ for $j>1$. These iterated brackets between the generators of $D_1$ induce a filtration:
 \[
D_1\subset D_2\subset D_3\subset\cdots
\]
 on $TM$. By definition  (\cite{Alekseevsky-Spiro-2001}), our distribution $D_1$ is {\sl regular} if for all $j \geqslant 1$, dimensions of the vector spaces $D_j(p)$ are independent of the choice of $p\in M$ and, moreover, there exists some (minimal) constant integer $\rho$ satisfying $D_\rho(p)=T_p M$ for each arbitrary point $p$.  In this case, $\rho$ is called the {\sl degree of nonholonomy} of the distribution $D_1$.

\begin{Definition} ({\it cf.} \cite[Definition 1.1]{Bel-Conj}).
\label{def-totally-nondegenerate}
An arbitrary (local) real analytic generic CR manifold $M$ is {\sl totally nondegenerate} (or {\sl maximally minimal}) of {\sl length} $\rho$ if the distribution $D_1=T^c M$ is regular with the {\it minimum possible degree} of nonholonomy $\rho$. The latter is equivalent to requiring that the length $\rho$ of the induced filtration:
\begin{equation}
\label{filtration}
T^cM=D_1\varsubsetneq D_2\varsubsetneq\cdots\varsubsetneq D_\rho=TM,
\end{equation}
be {\it minimum}.
\end{Definition}

In this definition, the term {\it "minimum possible length"} may seem somehow vague. In order to gain a deeper understanding, it may be helpful to recall some auxiliary algebraic notions.

\begin{Definition}
\label{def-fundamental-free-CR}
A finite dimensional graded Lie algebra:
\[
\frak g:=\frak g_{-\mu}\oplus\frak g_{-\mu+1}\oplus\cdots\oplus\frak g_{-1}
\]
is {\sl fundamental} if it can be generated by its subspace $\frak g_{-1}$, {\it i.e.} $\frak g_{-t}=[\frak g_{-1}, \frak g_{-t+1}]$, for each $t=2, \ldots, \mu$. In this case, two integers $\mu$ and $r:={\rm dim}\,\frak g_{-1}$ are called the {\sl length} and {\sl rank} of $\frak g$, respectively. Such a fundamental algebra $\frak g$ is:
\begin{itemize}
\item[(i)] {\sl free of length $\mu$} if on its iterated brackets of lengths $\leq\mu$ between the elements of $\frak g_{-1}$ there is no imposed relation except the {\it unavoidable} Jacobi identity and skew symmetry.
\item[(ii)] {\sl CR} (or {\sl pseudocomplex}) if it is equipped with a complex structure $J:\frak g_{-1}\rightarrow\frak g_{-1}$, satisfying:
\begin{equation}
\label{CR}
[J\,{\sf x}, J\,{\sf y}]=[{\sf x}, {\sf y}], \ \ \ \ \ \textrm{ for each} \ \ \ {\sf x,y}\in\frak g_{-1}.
\end{equation}
\item[(iii)] {\sl free CR of length $\mu$} if it is a quotient algebra of the form $\frak f/\frak i$ where $\frak f:=\frak f_{-\mu}\oplus\ldots\oplus\frak f_{-1}$ is a length $\mu$ free algebra endowed with a complex structure $J:\frak f_{-1}\rightarrow\frak f_{-1}$ and $\frak i$ is the ideal generated by all elements:
 \[
 [J\,{\sf x}, J\,{\sf y}]-[{\sf x}, {\sf y}] \ \ \ \ \ \  {\sf x,y}\in\frak f_{-1}.
 \]
\end{itemize}
\end{Definition}

When $\frak g=\frak g_{-\mu}\oplus\frak g_{-\mu+1}\oplus\cdots\oplus\frak g_{-1}$ is free, then each subspace $\frak g_{-t}$ has the maximum possible dimension as the $t$-th component of a fundamental algebra of the same rank. Furthermore, a free CR algebra is equivalently a fundamental algebra with no imposed relation on its iterated brackets except Jacobi identity, skew symmetry and the CR condition \thetag{\ref{CR}}.

 \medskip
Now, let us return to Definition \ref{def-totally-nondegenerate}. Consider a generic manifold $M$ of CR dimension $n$, codimension $k$ and with the associated regular distribution $T^cM$ with the degree of nonholonomy $\rho$. In agreement with the notations introduced above, we set $\frak D_{-1}=D_1$ and $\frak D_{-\ell}:= D_{\ell} / D_{\ell-1}$ for $\ell>1$. Equipping the graded space:
 \[
 \frak s(M):=\frak D_{-\rho}\oplus\frak D_{-\rho+1}\oplus\cdots\oplus\frak D_{-1}
 \]
 with the canonical bracket of vector fields converts it to a $(2n+k)$-dimensional fundamental algebra of length $\rho$ and rank $2n$. We call $\frak s(M)$ by the {\sl symbol algebra} of $M$. This Lie algebra is actually CR with respect to the complex structure defined on $\frak D_{-1}=T^cM$. Now, seeking the {\it minimum} length $\rho$ for the corresponding filtration \thetag{\ref{filtration}} forces each subspace $D_{\ell}$, and equivalently $\frak D_{-\ell}$, to have the {\it maximum possible dimension}. This together with the CR condition \thetag{\ref{CR}} are the only constraints imposed on the structure of $\frak s(M)$ by the requirements of Definition \ref{def-totally-nondegenerate}. Thus, we have the following characterization of total nondegeneracy.

\begin{Proposition}
\label{prop-prel}
A (local) real analytic generic CR manifold $M$ of CR dimension $n$ and codimension $k$ is {\sl totally nondegenerate} of {\sl length} $\rho$ if and only if its CR structure $T^cM$ is regular with the degree of nonholonomy $\rho$ and, moreover, the truncated graded algebra:
\[
\frak s(M) / \frak D_{-\rho}=\frak D_{-\rho+1}\oplus\cdots\oplus\frak D_{-1}
\]
is free CR of length $\rho-1$.
\end{Proposition}

\subsection{Levi-Tanaka prolongation}
\label{subsec-Tanaka}

The Tanaka theory of prolongations is a parallel algebraic version of Cartan's geometric approach to equivalence problems. The concept of transitivity is a cornerstone of this theory.

\begin{Definition}
\label{def-transitive}
Let $\frak g=\bigoplus_{i\in\mathbb Z}\frak g_i$ be a graded algebra and set $\frak g_-:=\bigoplus_{i\leqslant -1}\frak g_i$. The Lie algebra $\frak g$ is called {\sl transitive} whenever $[{\sf x}_i, \frak g_-]\not\equiv 0$ for each nonzero element ${\sf x}_i\in\frak g_i$ with $i\geqslant 0$.
\end{Definition}

Consider a fundamental algebra $\frak m:=\bigoplus_{-\mu\leqslant i\leqslant -1}\frak m_{i}$. Assume that it is CR provided by a certain complex structure $J$ on $\frak m_{-1}$. Define $\frak g^0(\frak m)$ as the algebra of all degree preserving derivations ${\sf d}:\frak m\rightarrow\frak m$ that respect $J$, {\it i.e.} ${\sf d}(J{\sf x})=J({\sf d}({\sf x}))$ for all ${\sf x}\in\frak m_{-1}$. Tanaka in \cite{Tanaka-main-1970} showed that there exists a unique transitive algebra $\frak g(\frak m):=\bigoplus_{-\mu\leq i}\frak g^i(\frak m)$ which is maximal among all transitive prolongations of $\frak m\oplus\frak g^0(\frak m)$. He also established a practical way to construct it. Accordingly, we set first $\frak g^{i}(\frak m):=\frak m_i$ for each $i\leq -1$. Assuming that the components $\frak g^{l'}(\frak m)$ are already constructed for any $l' \leqslant l - 1$, the $l$-th component $\frak g^l(\frak m)$ of the prolongation consists of $l$-shifted graded linear morphisms $\mathfrak m \to \mathfrak m \oplus \frak g^0(\frak m) \oplus \frak g^1(\frak m) \oplus \cdots\oplus \frak g^{ l-1}(\frak m)$ which are derivations:
\begin{equation*}
\frak g^l(\frak m)
=
\Big\{ {\sf d} \in
\bigoplus_{k\leqslant-1}\, {\sf Lin}\big(\frak g^k(\frak m),\,\frak g^{k+l}(\frak m)\big)
\colon
{\sf d}([{\sf y},\,{\sf z}]) =
[{\sf d}({\sf
y}),\,{\sf z}] + [{\sf y},\,{\sf d}({\sf z})],
\ \ \ \ \
\forall\, {\sf y},\,{\sf
z}\in\mathfrak m \Big\}.
\end{equation*}
Continuing in this fashion, we get the sought transitive prolongation:
\[
\frak g(\frak m):=\frak m\oplus\frak g^0(\frak m)\oplus\frak g^1(\frak m)\oplus\cdots.
\]
When the fundamental algebra $\frak m$ is {\sl nondegenerate}, that is: $[{\sf x}, \frak m_{-1}]\not\equiv 0$ for each nonzero element ${\sf x}\in\frak m_{-1}$, then the above process of construction will finally terminate (\cite{Tanaka-main-1970, Medori-Nacinovich-1997}).
In the case that $\frak m:=\frak s(M)$ is the symbol algebra of a certain CR manifold $ M$, then $\frak g(\frak m)$ is called the {\sl Levi-Tanaka algebra of $ M$} (\cite{Medori-Nacinovich-1997}). One notices that if $M$ is nondegenerate then so is its associated symbol algebra and hence the Levi-Tanaka algebra of $M$ is finite dimensional, as well.

\subsection{Totally nondegenerate models and the maximum conjecture}

 In the seminal work \cite{Beloshapka2004}, Beloshapka showed that after appropriate weight assignment to complex coordinates $(z_1, \ldots, z_n, w_1, \ldots, w_k)$, every totally nondegenerate submanifold of $\mathbb C^{n+k}$ of CR dimension $n$ and codimension $k$ can be represented locally around the origin as the graph of some $k$ real analytic defining functions ({\it cf.} \cite[Theorem 1.1]{Bel-Conj}):
 \begin{equation}
 \label{totally-nondegenerate}
\aligned
\begin{array}{l}
w_j=\Phi_j({\bf z},\overline {\bf z},\overline {\bf w})+{\rm O}([w_j]), \ \ \ \ \ \ \ \ \  j=1, \ldots, k
\end{array}
\endaligned
\end{equation}
where $[w_j]$ is the assigned weight number to $w_j$, where ${\rm O}(t)$ denotes some certain sum of monomials of weights $> t$ and where $\Phi_j$ is a weighted homogeneous complex-valued polynomial in terms of the coordinates ${\bf z}, \overline {\bf z}, w_j$ and other complex variables $w_\bullet$ of weights $[w_\bullet]<[w_j]$ (throughout this paper, we denote by $\bf z$ the $n$-tuple $(z_1, \ldots, z_n)$ and we use the similar notation for the variables $w$). Accordingly, Beloshapka introduced the specific totally nondegenerate CR manifold, represented by the corresponding defining polynomials:
  \begin{equation*}
\aligned
\begin{array}{l}
w_j=\Phi_j({\bf z},\overline {\bf z},\overline {\bf w}), \ \ \ \ \ \ \ \ \ \  j=1, \ldots, k
\end{array}
\endaligned
\end{equation*}
as the {\it model} of all  totally nondegenerate manifolds, represented as \thetag{\ref{totally-nondegenerate}}. He also established a practical way to construct the associated defining polynomials $\Phi_\bullet$ ({\it see} $\S$\ref{sec-Beloshapka-models}, below). We refer the reader to \cite[Theorem 14]{Beloshapka2004} for various {\it nice} features of these models.

Let $M\subset\mathbb C^{n+k}$, defined as the graph of $k$ polynomial functions $\Xi_j:=w_j-\Phi_j({\bf z}, {\bf \overline z}, {\bf \overline w})\equiv 0$ for $j=1,\ldots, k$, be a certain length $\rho$ Beloshapka's model with local coordinates $(z_1, \ldots, z_n,w_1,\ldots,w_k)$. A holomorphic vector field:
\[
{\sf X}:=\sum_{j=1}^n\,Z^j({\bf z},{\bf w})\frac{\partial}{\partial z_j}+\sum_{l=1}^k\,W^l({\bf z},{\bf w})\frac{\partial}{\partial w_l}
\]
is an {\sl infinitesimal CR automorphism} of $M$ if its real part is tangent to it, {\it i.e.} $({\sf X}+\overline{\sf X})|_{\Xi_j}\equiv 0$, $j=1,\ldots, k$. The collection of all such holomorphic vector fields forms a Lie algebra, denoted by $\frak{aut}_{CR}(M)$. As is known ({\it cf.} \cite{Beloshapka2004, Filomat}), this algebra is finite dimensional, of polynomial type and graded of the form:
\begin{equation}
\label{aut}
\frak{aut}_{CR}(M):=\underbrace{\frak g_{-\rho}\oplus\cdots\oplus\frak g_{-1}}_{\frak g_-}\oplus\,\frak g_{0}\oplus\underbrace{\frak g_{1}\oplus\cdots\oplus\frak g_{\mu}}_{\frak g_+}, \ \ \ \ \ \ \mu\,\in\,\mathbb N.
\end{equation}
The simply connected local Lie group ${\sf Aut}_{CR}(M)$, associated with $\frak{aut}_{CR}(M)$, comprises all CR-diffeomorphisms (or biholomorphisms), mapping $M$ to itself. Corresponding to \thetag{\ref{aut}}, one may write:
\begin{equation*}
\label{Aut}
{\sf Aut}_{CR}(M):={\sf G}_-\cdot {\sf G_0}\cdot {\sf G}_+.
\end{equation*}
As Beloshapka discovered (\cite[Proposition 3]{Beloshapka2004}), the Lie subgroup ${\sf G}_-$, associated with the subalgebra $\frak g_-$ of $\frak{aut}_{CR}(M)$, acts on $M$ freely and transitively and can be identified with $M$ through a certain diffeomorphism. The subgroup ${\sf G}_0$, associated with $\frak g_0$, comprises all {\it linear} automorphisms in the isotropy subgroup of ${\sf Aut}_{CR}(M)$ at the origin  while ${\sf G}_+$, associated with $\frak g_+$, comprises all {\it nonlinear} ones.

As is known ({\it cf.} \cite[Lemma 3.1]{Filomat}), the number $\rho$ in the index of the minimum component $\frak g_{-\rho}$ in \thetag{\ref{aut}} is actually equal to the length of $M$. But unfortunately, determining the value of the {\it maximum} index $\mu$ in this gradation is still an open question. Beloshapka in \cite{Beloshapka} conjectured that;

\begin{Conjecture}
\label{conjecture}
{\bf [Beloshapka's Maximum Conjecture]} Every totally nondegenerate CR model $M$ of length $\rho\geq 3$ has {\sl rigidity}; that is: in the gradation \eqref{aut} of $\frak{aut}_{CR}(M)$, the subalgebra $\frak g_+$ is trivial.
\end{Conjecture}

If this conjecture holds true, then Beloshapka's models of length $\rho\geq 3$ never have any {\it nonlinear} transformation in their CR automorphism groups. In the recent work \cite{Bel-Conj}, the author proved the maximum conjecture in CR dimension one by applying the celebrated classical approach of \'{E}lie Cartan toward solving (biholomorphic) equivalence problems. Moreover, Gammel and Kossovskiy in \cite{Gammel-Kossovskiy-2006} confirmed it in length $\rho=3$. Kossovskiy also proved this conjecture in \cite{Kossovskiy-2007} for length four model CR manifolds with {\it reflection}. Both of Kossovskiy's proofs are based on the structure of the {\it envelopes of holomorphy} of the corresponding models.

The main goal of this paper is to prove the maximum conjecture in two specific, but interesting, cases by means of some results and techniques arisen in the theory of (Levi-)Tanaka prolongations. First, we prove it in Section \ref{sec-rigidity} for the case of the so-called {\it full-models} which are those Beloshapka's models  that their associated symbol algebras are free CR. The result guarantees to have at least one totally nondegenerate model with rigidity in each arbitrary CR dimension $n$ and length $\rho\geq 3$. Next in Section \ref{sec-CR-dim-one}, we prove the conjecture in CR dimension one. As we already mentioned, the maximum conjecture in this CR dimension is actually proved before in \cite{Bel-Conj} by means of Cartan's {\it geometric} approach for solving equivalence problems. Nevertheless, here we provide a much shorter proof using the parallel {\it algebraic} approach of Tanaka.

Before concluding this preliminary section, it may be worth to notice\,\,---\,\,based on the extant experiences\,\,---\,\,that among the three sofar applied approaches, namely the envelope of holomorphy  (\cite{Gammel-Kossovskiy-2006, Kossovskiy-2007}), the Cartan theory (\cite{Bel-Conj}) and the Tanaka theory of prolongations (the current study), it seems that due to its convenience in comparison to other ones, {\it the algebraic approach of Tanaka is the best weapon to attack the maximum conjecture}. The results obtained in this paper are encouraging enough to merit further investigation by means of this approach.

\section{Beloshapka's totally nondegenerate models}
\label{sec-Beloshapka-models}

We review in this section the method of constructing defining equations of  Beloshapka's models.  Throughout the process of construction and to each complex coordinate ${\sf x}=z_i, w_j$, it will be assigned a {\sl weight} number, denoted by $[{\sf x}]$. We assign the same weight $[{\sf x}]$ to the conjugation $\overline{\sf x}$ and real and imaginary parts of $\sf x$, as well. Also to a certain monomial ${\sf x}_1^{\alpha_1}\cdots {\sf x}_n^{\alpha_n}$, the associated weight is defined as $\sum_{i=1}^n\,\alpha_i\,[{\sf x}_i]$. A polynomial is called {\sl weighted homogeneous} whenever all monomials appearing in it are of the same weight.

\begin{Definition} ({\it cf.} \cite{Krantz})
An arbitrary $\mathcal C^2$ complex function $f:\Omega\subset\mathbb C^n\rightarrow\mathbb C$ in terms of the coordinates $({\sf z}_1,\ldots, {\sf z}_n)$ is called {\sl pluriharmonic} on its domain if $\frac{\partial^2 f}{\partial {\sf z}_i\,\partial\overline{{\sf z}}_j}\equiv 0$ for each $i,j=1,\ldots, n$.
\end{Definition}

At the onset of constructing the sought defining polynomials of a model $M\subset\mathbb C^{n+k}$, of CR dimension $n$ and codimension $k$, we assign to all CR coordinates $z_i, i=1, \ldots, n$ the same weight $[z_i]=1$. The weights of the other complex variables $w_1, w_2, \ldots$\,\,---\,\,which are absolutely bigger than $1$\,\,---\,\,will be determined as follows, step by step.

At the first step, let $\mathcal N_2$ be a basis for the space of all non-pluriharmonic real-valued polynomials of the homogeneous weight $2$ in terms of the complex variables $z_i$ and their conjugations. A careful inspection shows that $\mathcal N_2$ comprises the terms:
\[
\mathcal N_2:=\{{\rm Re}\,(z_i\,\overline z_j), \ \  {\rm Im}\,(z_i\,\overline z_j), \ \ z_i\overline z_i, \ \ i, j=1, \ldots, n, \ \ i\neq j\}.
\]
 As the cardinality of this set is $k_2=n^2$, then we assign accordingly the weight $2$ to the first $k_2$ complex variables $w_1, \ldots, w_{k_2}$.

Let us label $k_2$ elements of $\mathcal N_2$ by ${\sf t}^2_1,\ldots,{\sf t}^2_{k_2}$. Define the next collection $\mathcal N_3$ as a basis for the space of all real-valued polynomials of the weight $3$, in terms of the variables $z_i, \overline z_i$, $i=1, \ldots, n$ and ${\rm Re}\,w_j$, $j=1, \ldots, k_2$, and non-pluriharmonic on the manifold represented as the graph of $k_2$ weighted homogeneous polynomials ${\rm Im}\,w_j= {\sf t}_j$, $j=1, \ldots, k_2$.
Again, a careful inspection shows that ({\it cf.} \cite{Gammel-Kossovskiy-2006}):
\[
\mathcal N_3:=\{{\rm Re}\,(z_i\, z_j\,\overline z_r), \ \ \ \ {\rm Im}\,(z_i\, z_j\,\overline z_r), \ \ \ i,j,r=1, \ldots, n\}.
\]
The cardinality of this basis is $k_3=n^3+n^2$, hence we assign the weight $3$ to the next $k_3$ complex variables $w_{k_2+1}, \ldots, w_{k_2+k_3}$.

Inductively, assume that $\mathcal N_{j_0}$ is the last constructed basis for some integer $j_0\in\mathbb N$. This means that at the moment all the complex variables $z_1, \ldots, z_n, w_1, \ldots, w_{\sf r}$ with ${\sf r}:=\sum_{i=2}^{j_0}\, k_i$ and $k_i:={\rm card}(\mathcal N_i)$ have received their weight numbers. Constructing the next basis $\mathcal N_{j_0+1}$, let us label $k_i$ basis elements of each $\mathcal N_i$ by $\{{\sf t}^i_1, {\sf t}^i_2, \ldots,{\sf t}^i_{k_i}\}$.  Also for each $\ell=2,\ldots, j_0$, let ${\bf w}_\ell=(w_{l}, \ldots, w_{l+k_\ell-1})$ be the $k_\ell$-tuple of the complex variables with the same weight $\ell$. Then, the sought collection $\mathcal N_{j_0+1}$ is defined as a basis for the space of all real-valued polynomials of weight $j_0+1$, in terms of the already weight determined variables $z_i, \overline z_j, {\rm Re}\,w_1, {\rm Re}\,w_2, \ldots, {\rm Re}\,w_{\sf r}$, $i,j =1, \ldots, n$ and non-pluriharmonic on the manifold represented as the graph of $\sf r$ weighted homogeneous polynomials:
\[
\aligned
({\rm Im}\,{\bf w}_j)^t
&=
\big(\, {\sf t}^j_1 \, \ldots \, {\sf t}^j_{k_j}\,\big )^t \ \ \ \ \ \ \ \ \ j=2,\ldots,j_0.
\endaligned
\]
 Here, ${\rm Im}\,{\bf w}_\ell$ is the $k_\ell$-tuple of the imaginary parts of ${\bf w}_\ell$.
If the cardinality of $\mathcal N_{j_0+1}$ is $k_{j_0+1}$ and proceeding as before, we assign the weight $j_0+1$ to all the next complex variables $w_{{\sf r}+1},\ldots,w_{{\sf r}+k_{j_0}}$.

 \subsection{Constructing the defining equations}

After assigning appropriate weights to complex variables, we are ready to construct defining polynomials of a Beloshapka's model $M\subset\mathbb C^{n+k}$ of CR dimension $n$ and codimension $k$. In this case, we need only the assigned weights to the complex coordinates $(z_1, \ldots, z_n,w_1,\ldots,w_k)$ of $M$. Thus, we only need to find the above bases $\mathcal N_i$ until we arrive at the stage $i=\rho$ as the smallest integer satisfying:
\begin{equation*}
 k\leqslant  k_2+\cdots+k_{\rho-1}+k_\rho.
 \end{equation*}
It is clear from the above process that the last coordinate variable $w_k$ of $M$ receives the maximum weight $\rho$. This integer is actually the length of the totally nondegenerate model $M$, in question.

Now, for each $\ell=2,\ldots, \rho-1$, consider the $k_\ell$-tuple ${\bf w}_\ell$ as above and let:
\begin{equation*}
A_\ell=\left(
         \begin{array}{ccc}
           a^\ell_{11} & \ldots & a^\ell_{1k_\ell} \\
           \vdots & \vdots & \vdots \\
            a^\ell_{k_\ell 1} & \ldots & a^\ell_{k_\ell k_\ell} \\
         \end{array}
       \right)
\end{equation*}
be some certain real $k_\ell\times k_\ell$ matrix of the maximum rank $k_\ell$. For $\ell=\rho$ and since the number of the present weight $\rho$ coordinate variables $w_\bullet$ of $M$ is $m=k-\sum_{i=2}^{\rho-1}k_i\leqslant k_\rho$, then we consider the $m$-tuple ${\bf w}_\rho=(w_{k-m+1},\ldots, w_k)$. Also let:
\[
A_\rho=\left(
         \begin{array}{ccc}
           a^\rho_{11} & \ldots & a^\rho_{1 k_\rho} \\
           \vdots & \vdots & \vdots \\
            a^\rho_{m 1} & \ldots & a^\rho_{m k_\rho} \\
         \end{array}
       \right)
\]
be a certain real $m\times k_\rho$ matrix of the maximum possible rank $m$. Then, the desired $k$ defining equations of $M$ might be represented in the following matrix form:

\begin{equation}
\label{model-matrix}
\aligned
({\rm Im}\,{\bf w}_2)^t
&=
A_2
\cdot
\big (\,{\sf t}^2_1\,  \ldots \, {\sf t}^2_{k_2}\,\big)^t,
\\
& \ \  \vdots
\\
({\rm Im}\,{\bf w}_{\rho})^t
&=
A_{\rho}
\cdot
\big(\,{\sf t}^\rho_1\,  \ldots \, {\sf t}^{\rho}_{k_{\rho}}\,\big )^t.
\endaligned
\end{equation}
It is important to notice that every collection of the entries $a^\ell_{ij}$ in the above matrices $A_2, \ldots, A_\rho$ may produce biholomorphically a distinct model of CR dimension $n$ and codimension $k$. Thus, in each fixed CR dimension and codimension we may find a large number of holomorphically inequivalent CR models.

For our length $\rho$ CR model $M$, the specific case that the already introduced $m\times k_\rho$ matrix $A_\rho$ is square, namely when the CR codimension $k$ of $M$ has the maximum possible value $k=\sum_{i=2}^\rho k_i$, is of particular interest. Thanks to \cite[Proposition 4(d)]{Beloshapka2004}, such a CR model can be characterized by the fact that {\it all}\,\,---\,\,even the last\,\,---\,\,components of its associated graded algebra $\frak s(M)$ have the maximum possible dimension. This, according to Proposition \ref{prop-prel}, is equivalent to state that $\frak s(M)$ is itself free CR.

\begin{Definition}
\label{def-full-models}
We call a length $\rho$ Beloshapka's model by a {\it full-model} if its associated symbol algebra is free CR. We denote such a model by $\mathbb M_\rho$.
\end{Definition}

\section{The rigidity of full-models}
\label{sec-rigidity}

In this section, we investigate the maximum conjecture in the case of full-models.
Let $\mathbb M_\rho$ be a fixed full-model of  CR dimension $n$ and codimension $k$. For each $\ell=1, \ldots, \rho$, consider $D_\ell\subset T\mathbb M_\rho$ and $\frak D_{-\ell}=D_\ell/ D_{\ell-1}$ as in Section \ref{sec-prel}. By definition, here the associated symbol algebra:
\[
\frak s(\mathbb M_\rho)=\frak D_{-\rho}\oplus\cdots\oplus\frak D_{-2}\oplus\frak D_{-1}
\]
is free CR of length $\rho$ and of rank $2n$. Let $\frak D_{-1}=T^c\mathbb M_\rho$ be generated by the basis $\{{\sf x}_1,\ldots,{\sf x}_n,J{\sf x}_1,\ldots,J{\sf x}_n\}$, where $J$ is the associated complex structure map. Equip this basis with a total ordering like:
\[
{\sf x}_1\prec\ldots\prec {\sf x}_n\prec J{\sf x}_1\prec\ldots\prec J{\sf x}_n.
\]
As discussed in Section \ref{sec-prel}, the subspace $\frak D_{-2}$ of $\frak s(\mathbb M_\rho)$ is constructed by all length two iterated brackets in the rank $2n$ real free algebra $\frak f$ which enjoy the CR condition \thetag{\ref{CR}}. Thus, a basis for $\frak D_{-2}$ is the second Hall basis of $\frak f$ ({\it see} \cite[$\S$4.1]{Reutenauer} for definitions):
\begin{equation}
\label{Hall-basis}
\aligned
\left\{
\begin{array}{l}
(i): \ \ \ \  \big[{\sf x}_i, {\sf x}_j\big], \ \ \ i<j,
\\
\\
(ii):\ \ \ \big[{\sf x}_i, J{\sf x}_j\big], \ \ \  i,j=1,\ldots, n,
\\
\\
(iii): \ \ \big[J{\sf x}_i, J{\sf x}_j\big], \ \ \ i<j,
\end{array}
\right.
\endaligned
\end{equation}
modulo the CR condition \thetag{\ref{CR}}. Consequently, $\frak D_{-2}$ can be generated by $n^2$ basis elements of the form $(i)$ and $(ii)$ or of the form $(ii)$ and $(iii)$, where we impose the extra restriction $i\leqslant j$ on the elements of $(ii)$\,\,---\,\,notice that \eqref{CR} forces the relation $[{\sf x}_i, J{\sf x}_j]=[{\sf x}_j, J{\sf x}_i]$. The next components $\frak D_{-j}$ with $j>2$ are constructed similarly as some subspaces of the free Lie algebra $\frak f$, imposed by the extra CR condition \thetag{\ref{CR}}.

\begin{Proposition}
\label{prop-Levi-Tanaka-full-model}
If $\rho>3$, then the Levi-Tanaka prolongation of the CR algebra $\frak s(\mathbb M_\rho)$ is trivial, that is: it does not include any component of positive degree.
\end{Proposition}

\proof
The proof of this assertion is actually an adjustment of Warhurst's main proof in \cite{Warhurst-2007} about the triviality of Tanaka prolongation of free algebras. We present here a sketch of it. Let:
\[
\frak g(\frak s(\mathbb M_\rho))=\frak s(\mathbb M_\rho)\oplus\frak g^0\oplus\frak g^1\oplus\cdots
\]
is the Levi-Tanaka prolongation of $\frak s(\mathbb M_\rho)$. We shall prove that $\frak g^1\equiv 0$. For this purpose, it suffices to show that if $u\in\frak g^1$ then $u({\sf y})({\sf z})=0$ for each basis elements ${\sf y}, {\sf z}={\sf x}_1,\ldots,{\sf x}_n,J{\sf x}_1,\ldots,J{\sf x}_n$ of $\frak D_{-1}$ ({\it see} the paragraph after Theorem 1 in \cite{Warhurst-2007}). Seeking the simplicity, let us denote $u({\sf y})(\sf z)$ by $u{\sf y}\sf z$. Since for each basis element $\sf y$ of $\frak D_{-1}$, we have $u{\sf y}\in\frak g^0$ and, by definition, all the elements of $\frak g^0$ respect the associated structure map $J$, then it is sufficient to show:
\[
u\,{\sf x}_i\,{\sf x}_j=0 \ \ \ \ {\rm and} \ \ \ \ u\,J{\sf x}_i\,{\sf x}_j=0
\]
for each $i,j=1,\ldots,n$. We prove here only the former equality. Proof of the latter one can be obtained by a similar argument. For two fixed indices $i,j=1,\ldots,n$, assume that:
\begin{equation}
\label{uxx}
\aligned
u\,{\sf x}_i\,{\sf x}_i&=a_1\,{\sf x}_i+b_1\,{\sf x}_j+{\bf c}_1\,{\bf x}+{\bf d}_1\,J{\bf x},
\\
u\,{\sf x}_i\,{\sf x}_j&=a_2\,{\sf x}_i+b_2\,{\sf x}_j+{\bf c}_2\,{\bf x}+{\bf d}_2\,J{\bf x},
\\
u\,{\sf x}_j\,{\sf x}_i&=a_3\,{\sf x}_i+b_3\,{\sf x}_j+{\bf c}_3\,{\bf x}+{\bf d}_3\,J{\bf x},
\\
u\,{\sf x}_j\,{\sf x}_j&=a_4\,{\sf x}_i+b_4\,{\sf x}_j+{\bf c}_4\,{\bf x}+{\bf d}_4\,J{\bf x}.
\endaligned
\end{equation}
For the sake of brevity, we denote simply by ${\bf c}_\bullet{\bf x}$ and ${\bf d}_\bullet J{\bf x}$ some two certain combinations $\sum_{r\neq i, j}c_{\bullet r} {\sf x}_r$ and $\sum_{r'\neq i,j}d_{\bullet r'} J{\sf x}_{r'}$. Applying Lemma 1 of \cite{Warhurst-2007}, for $r=\rho$, $X={\sf x}_i$ and $Y={\sf x}_j$ gives:
\begin{equation}
\label{eq1}
\aligned\footnotesize
0=u({\rm ad}^\rho_{{\sf x}_i}\,{\sf x}_j)&=\big(\rho\,b_2-b_3+a_1\,C(\rho,0)\big)\,{\rm ad}^{\rho-1}_{{\sf x}_i}\,{\sf x}_j
\\
&+b_1\,\sum_{j=0}^{[\frac{\rho-1}{2}]}\,\big(C(\rho,j)-C(\rho,\rho-2-j)\big)\,\big[{\rm ad}^j_{{\sf x}_i}\,{\sf x}_j, {\rm ad}^{\rho-2-j}_{{\sf x}_i}\,{\sf x}_j\big]
\\
& +\big(\rho\,{\bf c}_2-{\bf c}_3\big)\,{\rm ad}^{\rho-1}_{{\sf x}_i}\,{\bf x}+\big(\rho\,{\bf d}_2-{\bf d}_3\big)\,{\rm ad}^{\rho-1}_{{\sf x}_i}\,J{\bf x}
\\
&+{\bf c}_1\,\sum_{j=0}^{[\frac{\rho-1}{2}]}\,\big(C(\rho,j)-C(\rho,\rho-2-j)\big)\,\big[{\rm ad}^j_{{\sf x}_i}\,{\bf x}, {\rm ad}^{\rho-2-j}_{{\sf x}_i}\,{\sf x}_j\big]
\\
&+{\bf d}_1\,\sum_{j=0}^{[\frac{\rho-1}{2}]}\,\big(C(\rho,j)-C(\rho,\rho-2-j)\big)\,\big[{\rm ad}^j_{{\sf x}_i}\,J{\bf x}, {\rm ad}^{\rho-2-j}_{{\sf x}_i}\,{\sf x}_j\big],
\endaligned
\end{equation}
where:
\[
C(r,j)=(j+1)\,
\left(
  \begin{array}{c}
    r \\
    j+2 \\
  \end{array}
\right).
\]
According to the discussion before this proposition, all appearing Lie brackets in the above expression can be included in a basis of $\frak s(\mathbb M_\rho)$. Thus, from the second, fourth and fifth lines of this equality we have:
\[
b_1={\bf c}_1={\bf d}_1=0.
\]
Similarly, by considering the equality $0=u({\rm ad}^\rho_{{\sf x}_j}\,{\sf x}_i)$, one obtains:
\[
a_4={\bf c}_4={\bf d}_4=0.
\]
Moreover, we have:
\begin{equation}
\label{eq-ab}
\rho\,b_2-b_3+a_1\,C(\rho,0)=0 \ \ \ \ {\rm and} \ \ \ \rho\,b_3-b_2+b_4\,C(\rho,0)=0,
\end{equation}
where the first equality is obtained from the first line of \thetag{\ref{eq1}} and the second one can be found by a similar argument for $0=u({\rm ad}^\rho_{{\sf x}_j}{\sf x}_i)$. So far, we have found the values of the six, from sixteen unknowns in \thetag{\ref{uxx}}. In order to fined the the remained ones, let us consider the equality $0=u({\rm ad}^{\rho-2}_{{\sf x}_i}[{\sf x}_j, [{\sf x}_i, {\sf x}_j]])$. For this purpose, applying once again Lemma 1 of \cite{Warhurst-2007} gives that:
\[
\aligned
u({\rm ad}^{\rho-2}_{{\sf x}_i}\,[{\sf x}_j, [{\sf x}_i, {\sf x}_j]])&={\rm ad}^{\rho-2}_{{\sf x}_i}\,u([{\sf x}_j, [{\sf x}_i, {\sf x}_j]])+(s-2)\,{\rm ad}^{\rho-3}_{{\sf x}_i}\,u({\sf x}_i)[{\sf x}_j, [{\sf x}_i, {\sf x}_j]]
\\
&+\sum_{j=0}^{\rho-4}\,C(\rho-2,j)\,\big[{\rm ad}^{j}_{{\sf x}_i}\,u\,{\sf x}_i\,{\sf x}_i,{\rm ad}^{\rho-4-j}_{{\sf x}_i}[{\sf x}_j, [{\sf x}_i, {\sf x}_j]]\big].
\endaligned
\]
Performing (somehow) long computations like those at the lower half of the page 65 of \cite{Warhurst-2007}, one obtains:
\begin{equation}
\label{eq2}
\aligned
0=u({\rm ad}^{\rho-2}_{{\sf x}_i}&\,[{\sf x}_j, [{\sf x}_i, {\sf x}_j]])=\big(2\,a_3+b_4+(\rho-3)\,a_2\big)\,{\rm ad}^{\rho-1}_{{\sf x}_i}\,{\sf x}_j
\\
&+(\rho-2)\big(2\,b_2+\frac{\rho-1}{2}a_1\big)\,{\rm ad}^{\rho-3}_{{\sf x}_i}\,\big[{\sf x}_j, [{\sf x}_i, {\sf x}_j]\big]
\\
&+\big(2\,{\bf c}_3+{\bf c}_2\,(\rho-3)\big)\,{\rm ad}^{\rho-3}_{{\sf x}_i}\,\big[{\bf x}, [{\sf x}_i, {\sf x}_j]\big]
+\big({\bf c}_2\,(\rho-1)-2\,{\bf c}_3\big)\,{\rm ad}^{\rho-3}_{{\sf x}_i}\,\big[{\sf x}_j, [{\sf x}_i, {\bf x}]\big]
\\
&+\big(2\,{\bf d}_3+{\bf d}_2\,(\rho-3)\big)\,{\rm ad}^{\rho-3}_{{\sf x}_i}\,\big[J{\bf x}, [{\sf x}_i, {\sf x}_j]\big]
+\big({\bf d}_2\,(\rho-1)-2\,{\bf d}_3\big)\,{\rm ad}^{\rho-3}_{{\sf x}_i}\,\big[{\sf x}_j, [{\sf x}_i, J{\bf x}]\big].
\endaligned
\end{equation}
By an appropriate refinement, if necessary, on the total ordering of the basis elements in $\frak D_{-1}$, all the appearing brackets in this expression can be included in some basis of $\frak D_{-\rho}$ and hence they are linearly independent. Then, equating to zero all the coefficients among the last two lines of \thetag{\ref{eq2}} immediately implies that:
\[
{\bf c}_2={\bf c}_3={\bf d}_2={\bf d}_3=0.
\]
Furthermore, from the first two lines of this expression and also by applying similar argument for $0=u({\rm ad}^{\rho-2}_{{\sf x}_j}[{\sf x}_i, [{\sf x}_j, {\sf x}_i]])$, we receive the following four equations:
\begin{equation}
\label{eq-ab-2}
\aligned
2\,a_3+b_4+(\rho-3)\,a_2=0, \ \ \ \ \ \ \ \ \ \ \ \ \ \ \ \ \ 2\,b_2+\frac{\rho-1}{2}a_1=0,
\\
2\,b_2+a_1+(\rho-3)\,b_3=0, \ \ \ \ \ \ \ \ \ \ \ \ \ \ \ \ \ 2\,a_3+\frac{\rho-1}{2}b_4=0.
\endaligned
\end{equation}
Let us display the six algebraic equations in \thetag{\ref{eq-ab}} and \thetag{\ref{eq-ab-2}} in the following matrix form:
\[
\aligned
\left(
  \begin{array}{cccccc}
    \frac{\rho\,(\rho-1)}{2} & 0 & \rho & 0 & -1 & 0 \\
    0 & 0 & -1 & 0 & \rho & \frac{\rho\,(\rho-1)}{2} \\
    0 & \rho-3 & 0 & 2 & 0 & 1 \\
    \frac{\rho-1}{2} & 0 & 2 & 0 & 0 & 0 \\
    1 & 0 & 2 & \rho-3 & 0 & 0 \\
    0 & 0 & 0 & 2 & 0 & \frac{\rho-1}{2} \\
  \end{array}
\right)
\cdot
\left(
  \begin{array}{c}
    a_1 \\
    a_2 \\
    b_2 \\
    a_3 \\
    b_3 \\
    b_4 \\
  \end{array}
\right)
=
0.
\endaligned
\]
Our computations show that the determinant of the above $6\times 6$ matrix is nonzero when $\rho>3$, as was assumed. Then, the only possible value for the remaining unknowns in \thetag{\ref{uxx}} is zero, as well. This completes the proof.
\endproof

Now, we are ready to present the main result of this section.

\begin{Theorem}
\label{cor-rigidity-free-case}
Let $\mathbb M_\rho$ be a full-model of length $\rho\geq 3$. Then, the gradation of the Lie algebra $\frak{aut}_{CR}(\mathbb M_\rho)$ of its infinitesimal CR automorphisms (cf. \eqref{aut}) does not include any positive homogeneous component. In other words, $\mathbb M_\rho$ has rigidity.
\end{Theorem}

\proof
For $\rho=3$, the assertion is proved in \cite{Gammel-Kossovskiy-2006}. Then, let us assume that $\rho>3$. Let ({\it cf.} \eqref{aut}):
\[
\frak{aut}_{CR}(\mathbb M_\rho)=\frak g_-\oplus\frak g_0\oplus\frak g_+.
\]
As stated in \cite[p. 483]{Beloshapka2004}, here we have ${\rm dim}\frak g_0=2n^2$ and thus:
 \begin{equation}
 \label{dim-aut}
 {\rm dim}\,\frak{aut}_{CR}(\mathbb M_\rho)={\rm dim}\,\frak g_-+2\,n^2+{\rm dim}\,\frak g_+.
 \end{equation}
On the other hand, the fundamental Lie subalgebra $\frak g_-$ can be identified with the symbol algebra $\frak s(\mathbb M_\rho)$ ({\it cf.} \cite[Proposition 4(c)]{Beloshapka2004}) and hence they have the same Levi-Tanaka transitive prolongation. Hence, according to Proposition \ref{prop-Levi-Tanaka-full-model}, this Levi-Tanaka prolongation of $\frak g_-$ is trivial of the form $\frak g(\frak g_-):=\frak g_-\oplus\frak g^0$. By definition, the maximum possible dimension for $\frak g^0$ is $2n^2$ and thus ${\rm dim}\frak g_-+2n^2$ is the maximum possible dimension for $\frak g(\frak g_-)$. Now, according to \cite[Proposition 10.6]{Tanaka-main-1970} we have:
 \begin{equation*}
{\rm dim}\,\frak{aut}_{CR}(\mathbb M_\rho)\leqslant {\rm dim}\,\frak g(\frak g_-)\leqslant{\rm dim}\frak g_-+2n^2.
 \end{equation*}
Now, comparing this inequality with \thetag{\ref{dim-aut}} immediately implies that $\frak g_+$ is zero dimensional and this completes the proof.
 \endproof

\begin{Remark}
At each arbitrary CR dimension $n$ and length $\rho\geq 3$, this result guarantees that there exists at least one CR model that enjoys the rigidity condition.
\end{Remark}

\section{The rigidity of models in CR dimension one}
\label{sec-CR-dim-one}

In this section, we attempt to prove the maximum conjecture in CR dimension one.
As mentioned, for an arbitrary totally nondegenerate model $M$ the Lie algebra of its infinitesimal CR automorphisms is graded of the form $\frak{aut}_{CR}(M):=\frak g_-\oplus\frak g_0\oplus\frak g_+$ and, moreover, the associated Lie group ${\sf Aut}_{CR}(M):={\sf G}_-\cdot {\sf G}_0\cdot{\sf G}_+$ of this algebra consists of all automorphisms $h:M\rightarrow M$ satisfying $h_\ast(T^cM)=T^cM$. According to \cite[Proposition 3]{Beloshapka2004}, our CR model $M$ is an ${\sf Aut}_{CR}(M)$-homogeneous space and there exists a certain diffeomorphism:
\begin{equation}
\label{diff}
\Gamma: M \longrightarrow \frac{{\sf Aut}_{CR}(M)}{{\sf G_0}\cdot\sf G_+}\cong{\sf G_-} \ \ \ \ \ \ \ {\rm with} \ \  \ \ \Gamma(0)={\tt e},
\end{equation}
where $\tt e$ is the identity element of $\sf G_-$.

Since the subalgebra $\frak g_-$ can be identified with the symbol algebra $\frak s(M)$ of $M$ ({\it cf.} \cite[Proposition 4(c)]{Beloshapka2004}) then it is CR with respect to the complex structure map $J$ defined on $T^cM$. Let us denote by ${\sf Aut}_J(\frak g_-)$ the Lie group of all automorphisms of $\frak g_-$, preserving the gradation and respecting the complex structure map $J$. As mentioned at the page 72 of \cite{Tanaka-main-1970} ({\it see also} \cite[Proposition 1.120]{Knap}), the associated Lie algebra  to this group is actually the zero component $\frak g^0(\frak g_-)$ of the Levi-Tanaka algebra of $M$. On the other hand, in this case that $\sf G_-$ is connected and simply connected, two automorphism Lie groups ${\sf Aut}({\sf G_-})$ and ${\sf Aut}(\frak g_-)$ are isomorphic through the map:
\[
\aligned
\Phi: {\sf Aut}({\sf G_-})&\longrightarrow {\sf Aut}(\frak g_-)
\\
f&\mapsto f_{\ast {\tt e}},
\endaligned
\]
where $f_{\ast \tt e}$ is the differentiation of $f:{\sf G_-}\rightarrow {\sf G_-}$ at the identity element $\tt e$ of ${\sf G}_-$ ({\it cf.} \cite{Hochschild-1952}). Define ${\sf Aut}_J({\sf G_-})\subset {\sf Aut}({\sf G_-})$ as the subgroup of all automorphisms $f$ with $f_{\ast {\tt e}}\in{\sf Aut}_J(\frak g_-)$. Then, clearly we have:

\begin{Lemma}
Through the above isomorphism $\Phi$, two Lie groups ${\sf Aut}_J(\frak g_-)$ and ${\sf Aut}_J({\sf G_-})$ are isomorphic. As a result, the Lie algebra of ${\sf Aut}_J({\sf G_-})$ is $\frak g^0(\frak g_-)$.
\end{Lemma}

Granted the fact that ${\sf G_0}$ is a subgroup of ${\sf Aut}_{CR}(M)$ and since we can identify $M$ with ${\sf G_-}$ through the above diffeomorphism $\Gamma$, this lemma implies that the Lie algebra $\frak g_0$, associated with ${\sf G}_0$, can be regarded as a subalgebra of $\frak g^0(\frak g_-)$. Now, we are ready to prove the maximum conjecture in CR dimension one.

\begin{Theorem}
Every totally nondegenerate CR model $M$ of CR dimension one and length $\rho\geqslant 3$ has rigidity.
\end{Theorem}

\proof
According to \cite[$\S$5.6]{Medori-Nacinovich-1997}, the Levi-Tanaka prolongation of the length $\rho$ CR subalgebra $\frak g_-:=\frak g_{-\rho}\oplus\ldots\oplus\frak g_{-1}$ of $\frak{aut}_{CR}(M)$ with the 2-dimensional generator $\frak g_{-1}$ is of the simple form:
\[
\frak g(\frak g_-)=\frak g_-\oplus\frak g^0(\frak g_-).
\]
Since $\frak g_0$ is a subalgebra of $\frak g^0(\frak g_-)$ and $\frak{aut}_{CR}(M)$ is transitive then we can regard $\frak{aut}_{CR}(M)$ as a subalgebra of the Levi-Tanaka prolongation $\frak g(\frak g_-)$. This implies that as well as $\frak g(\frak g_-)$, it can not admit any component of positive degree in its gradation.
\endproof

\subsection*{Acknowledgment}
The author gratefully thanks the referee for the helpful suggestions concerning the presentation of this paper. The research of the author was supported in part by the grant from IPM, no. 96510425.

\end{document}